\documentclass[12pt]{article}
\usepackage{latexsym,amsmath,amssymb}
\begin{document}

\begin{center}
{\large \bf On the Robust Stability of Stationary Solutions 
\\ 
to a Class of Mathieu-Type Equations}
\\[5pt]
G.V. Demidenko${}^{1,2}$, K.S. Myagkikh${}^2$
\\[5pt]
${}^1$
Sobolev Institute of Mathematics SB RAS, Acad. Koptyug avenue, 4, 
\\
Novosibirsk, 630090, Russia
\\
${}^2$
Novosibirsk State University, Pirogov str., 1, Novosibirsk, 630090 Russia
\\
demidenk@math.nsc.ru; k.myagkikh@g.nsu.ru
\end{center}

\begin{quote}
We consider a class of nonlinear ordinary differential equations of the second order 
with parameters. We establish conditions for perturbations of the coefficients of the equation 
under which the zero solution is asymptotically stable.
Estimates for attraction sets of the zero solution 
and estimates of the stabilization rate of solutions at infinity are obtained.
Using these results, theorems on the robust stability of stationary solutions are proven.
\\[5pt]
\\[5pt]
\end{quote}

\vskip20pt

{\bf 1. Introduction}

\vskip10pt
	
In the paper we study the robust stability of stationary solutions 
to differential equations of the following form
\begin{equation} \label{1}
y'' + \alpha \mu y' + (\beta \mu^2 + \mu \varphi(t)) f(y) = 0, 
\end{equation}
where
$\alpha, \beta > 0$
are constant coefficients, 
$\varphi(t)$ is a continuous function such that
\begin{equation} \label{2}
\varphi(t+T) \equiv \varphi(t), \quad
\int\limits_0^T\varphi(s)ds = 0,
\end{equation} 
$\mu > 0$
is a small parameter.
	
The study of some Mathieu-type equations leads to the consideration of equations 
of the form~(\ref{1}). In particular, the study of motion of 
an inverted pendulum whose suspension point performs high-frequency 
harmonic oscillations can be reduced to the consideration of such equations.
Indeed, in this case, the equation of motion of the pendulum has the form \cite{Bo} 
\begin{equation} \label{3}
y'' + \lambda y' + \frac{g - a \omega^2 \sin\,(\omega\tau)}{l} \sin\,y = 0, 
\end{equation} 
where
$y(\tau)$
is the deviation angle of the pendulum from the lower vertical equilibrium position,
$l$
is the length of the pendulum,
$g$
is the gravitational acceleration,
$\lambda$
is the coefficient of friction,
$x(\tau) = a \sin\,(\omega \tau)$
is the motion equation of the suspension point, herewith
$\omega$
is the frequency of its oscillations and
$a$
is the oscillation amplitude.
Therefore, moving on to the ``fast time'' by changing the variable
$t = \omega \tau$  
and denoting
$$
\mu =  \frac{a}{l}, \qquad \beta = \frac{g l}{a^2 \omega^2}, \qquad \alpha = \frac{\lambda l}{a\omega}, 
$$
(\ref{3}) can be rewritten as follows
\begin{equation} \label{4}
y'' + \alpha \mu y' + (\beta \mu^2 - \mu \sin\,t) \sin\,y = 0.
\end{equation} 
Linearizing (\ref{4}) in the neighborhood of the solution
$y(t) \equiv \pi$, 
we obtain the equation
$$
\hat y'' + \alpha \mu \hat y' - (\beta \mu^2 - \mu \sin\,t) \hat y = 0
$$
which is reduced to the Mathieu equation
$$
v'' - \left(\beta \mu^2 + \frac{\alpha^2 \mu^2}{4} - \mu \sin\,t\right) v = 0,
$$
where 
$v(\tau) = e^{\alpha \mu \tau/2}\hat y(\tau)$,
	
Remind that the study of (\ref{3}) gives an understanding of the unexpected effect
of behavior of the inverted pendulum when the suspension point oscillates.
Namely, for sufficiently high oscillation frequency
$\omega \gg 1$
and sufficiently small amplitude of the oscillations of the suspension point,
the upper equilibrium position of the pendulum becomes stable.
This effect was first predicted back by A.~Stephenson in 1908 (see~\cite{St}).
A rigorous proof of this fact was obtained by N.N.~Bogolyubov in 1942 (see~\cite{Bo}).
He showed that the stationary solution
$y(\tau) \equiv \pi$  
corresponding to the upper vertical equili\-brium position of the pendulum 
is asymptotically stable for 
$$
a \omega > \sqrt{2 g l}, \qquad \frac{a}{l} \ll 1.
$$ 
An explanation of this effect from a physical point of view 
is contained in the works of P.L.~Kapitsa~\cite{Ka1, Ka2}.
		
At present, there are various approaches to the proof of N.N.~Bogolyubov's theorem.
The classic proof is carried out by using the averag\-ing method
and its modification~\cite{BoMi, Mi}.
In \cite{DeMa04, DeDuMa}, an approach proposed by G.V.~Demidenko and I.I.~Matveeva
for the study of the stability of solutions to differential equations 
is used to solve the stability problem of the inverted pendulum.
This approach is based on the use of the criterion of the asymptotic stability of solutions
to linear differential equations with periodic coefficients in terms of the solvability
of boundary value problems for the Lyapunov differential equation~\cite{DeMa01}.
This approach allows not only to establish the stability,
but also to estimate the attraction sets as well as the stabilization rate 
of the solutions as 
$t \to \infty$.
		
In this paper, using the approach~\cite{DeMa04} and the technique described in~\cite{DeDuMa},
we obtain results on the robust stability of stationary solutions to 
equations of the form~(\ref{1}).
Namely, we indicate estimates for perturbations of the coefficients
under which the stationary solutions to the equation
$$
\bar y'' + (\alpha + \Delta\alpha) \mu \bar y' + \left((\beta + \Delta\beta) \mu^2 
+ \mu (\varphi(t) + \Delta\varphi(t))\right) f(\bar y) = 0  
$$
are asymptotically stable for
$0 < \mu \ll 1$. 
In particular, from these estimates the robust stability of the upper equilibrium position 
of the pendulum whose suspension point makes high-frequency oscillations follow.	
		
\vskip20pt

{\bf 2. Preliminary information}

\vskip10pt
	
Studying the stability of stationary solutions to (\ref{1}),
we use a number of results from~\cite{DeMa04}, \cite{DeMa01}.
We give the necessary assertions below.
	
Consider the linear homogeneous system of differential equations with periodic coefficients
\begin{equation} \label{5}
\dfrac{dy}{dt} = A(t) y, 
\qquad t \in R,
\end{equation}
where
$A(t)$
is a continuous $T$-periodic $(n \times n)$-matrix.

At first we remind the spectral criterion of the asymptotic stability
of the zero solution to (\ref{5}) (see, for example, \cite{DaKr}). 

\vskip10pt
	
{\bf Theorem 1.} 
{\it The zero solution to (\ref{5}) is asymptotically stable
if and only if the spectrum of the monodromy matrix
$Y(T)$ 
lies inside the unit disk
$\gamma = \{ \lambda: |\lambda| < 1 \}$. 
}
	
Now we give the criterion of the asymptotic stability of the zero solution to (\ref{5}),
which is formulated in terms of the solvability of the following boundary 
value problem~\cite{DeMa01}
\begin{equation} \label{6}
	\begin{array}{l}
		\displaystyle
		\frac{dH}{dt} + HA(t) + A(t)^*H = -C(t), \quad 0 \le t \le T,\\
		H(0) = H(T).
	\end{array}
\end{equation}

\vskip10pt
	
{\bf Theorem 2.} 
{\it The zero solution to (\ref{5}) is asymptotically stable
if and only if the boundary value problem~(\ref{6}) with
$C(t) = C^*(t) > 0$, $t \in [0,T]$,
has a Hermitian solution
$H(t)$ 
such that
$H(0) > 0$. 
}
	
Note that if the zero solution to (\ref{5}) is asymptotically stable,
then there exists a unique Hermitian solution
$H(t)$ 
for any continuous matrix
$C(t)$, 
and an explicit formula can be written for $H(t)$ (see \cite{DeMa01}, \cite{De09}).
In particular, as shown in \cite{DeMa01},
$H(t) > 0$ 
on the whole segment
$[0,T]$.
	
We now extend the matrices
$C(t)$ and $H(t)$
$T$-periodically on the whole real axis. 
Then, keeping the same notation, 
$H(t)$
can be written as
\begin{equation} \label{7}
H(t) = \left(Y^*(t)\right)^{-1} 
\left(\int\limits^\infty_t Y^*(s)C(s)Y(s) ds\right) Y^{-1}(t), 
\end{equation}
where
$Y(t)$
is the matrizant for~(\ref{5}) (see \cite{De09}).
Obviously, this matrix is a $T$-periodic solution to the Lyapunov differential equation
$$
\frac{dH}{dt} + HA(t) + A(t)^*H = -C(t), \quad C(t+T) = C(t), \quad t \in R.  
$$
		
Using the matrix 
$H(t)$ 
defined by (\ref{7}), we can obtain an estimate of the stabilization rate 
of solutions to (\ref{5}) as
$t \to \infty$.
We now formulate the following theorem from~\cite{DeMa01}.

\vskip10pt
		
{\bf Theorem 3.} 
{\it Suppose that the zero solution to (\ref{5})
is asymptotically stable and $T$-periodic matrix
$H(t)$
is defined in~(\ref{7}).
Then, for solutions to (\ref{5}), the estimate holds
$$
\|y(t)\|^2 \le \frac{\|H(0)\|}{h_{\min}(t)} \|y(0)\|^2 
\exp\left(-\int\limits^t_0 \frac{c_{\min}(s)}{\|H(s)\|}ds \right), \quad t \ge 0, 
$$
where
$h_{\min}(t) > 0$ and $c_{\min}(t) > 0$
are minimal eigenvalues of the matrices
$H(t)$ and $C(t)$,
respectively.
}
	
This estimate is an analogue of the Krein estimate for solutions 
to systems of ordinary differential equations with constant coefficients.
	
\vskip20pt

{\bf 3. Linearized equations}

\vskip10pt
	
Consider the differential equation of the form (\ref{1}) 
$$
y'' + \alpha \mu y' + (\beta \mu^2 + \mu \varphi(t)) f(y) = 0, 
$$
where the function  
$\varphi(t)$ 
satisfies (\ref{2}), the function 
$f(y)$ 
is smooth and there exists a constant 
$\gamma$ 
such that 
\begin{equation} \label{8}
	f(\gamma) = 0, \qquad \frac{d}{dy}f|_{y = \gamma} < 0. 
\end{equation}
Obviously, the function 
$y(t) \equiv \gamma$ 
is a stationary solution to (\ref{1}).  
	
Consider the linearized equation for (\ref{1}) in a neighborhood of the solution 
$y(t) \equiv \gamma$. 
Introducing the notation  
$$
	\hat\beta = \beta \frac{d}{dy}f|_{y = \gamma}, 
	\qquad 
	\hat\varphi(t) = \varphi(t) \frac{d}{dy}f|_{y = \gamma}, 
$$ 
it can be written in the form 
\begin{equation} \label{9}
	y'' + \alpha \mu y' + (\hat\beta \mu^2 + \mu \hat\varphi(t)) y = 0.  
\end{equation} 
Since 
$\beta > 0$, 
then from (\ref{8}) it follows that  
$\hat\beta < 0$.  
Taking into account the equality 
$$
\int\limits^{T}_0 \hat\varphi(s) ds = 0, 
$$ 
we obtain the next theorem for (\ref{9}) from \cite{DeDuMa}.

\vskip10pt

{\bf Theorem 4.} 
{\it If 
$\alpha > 0$ and
\begin{equation} \label{10}
	\frac{1}{T} \int\limits^{T}_0 
	\left( \int\limits^{\tau}_0 \hat\varphi(s)\, ds \right)^2 d\tau
	> \left(\frac{1}{T} \int\limits^{T}_0 \tau \hat\varphi (\tau) d\tau \right)^2 - \hat\beta,
\end{equation} 
then the zero solution to (\ref{9}) is asymp\-totically stable for sufficiently small 
$\mu>0$.
}
	
As mentioned above, the equation of motion of the pendulum (\ref{3}) can be rewritten 
in the form (\ref{4}). 
Therefore, linearizing it in the neighborhood of the solution 
$y(t) = \pi$, 
we come to (\ref{9}) with 
$$
\hat\beta = -\beta = -\frac{gl}{a^2\omega^2}, \qquad \hat\varphi(t) = \sin\,t. 
$$ 
Calculating the integrals in~(\ref{10}), we obtain the inequality 
$$ 
a^2\omega^2 > 2 gl. 
$$  
Hence, by Theorem 4 the zero solution to (\ref{9}) is asymptotically stable
for sufficiently small 
$\mu>0$. 
The above inequality is N.N.~Bogolyubov's condition for the stability of oscillations 
of the inverted pendulum with the vibrating suspension point.  
	
Observe that Theorem 4 does not give a change range of the small parameter   
$\mu$ 
under which the zero solution to~(\ref{9}) is asymptotically stable. 
Later we will indicate the number 
$\mu_0 > 0$ 
such that the zero solution to~(\ref{9}) is asymptotically stable for  
$\mu \in (0, \mu_0]$.  

Let 
$y(t)$ 
be a solution to~(\ref{9}) and 
$$
v(t) =  
\begin{pmatrix}
	v_1(t) \\
	v_2(t)    
\end{pmatrix} 
= \begin{pmatrix}
	y(t) \\
	y'(t)    
\end{pmatrix}. \\ 
$$
Then (\ref{9}) is equivalent to the system  
\begin{equation} \label{11}
	\frac{d}{dt}v = A(t,\mu) v, 
\end{equation}
where 
$$
A(t,\mu) = 
\begin{pmatrix}
	0 & 1 \\
	-(\hat\beta \mu^2 + \mu \hat\varphi(t)) & -\alpha \mu 
\end{pmatrix}. 
$$ 
Following \cite{DeDuMa}, we change the variables
\begin{equation} \label{12}
	v_1(t) = (1 + \mu a(t)) u_1(t), \qquad v_2(t) = \mu b(t) u_1(t) + \mu c(t) u_2(t),
\end{equation}
where 
\begin{equation} \label{13}
	a(t) = -\int\limits^t_{t_2} \int\limits^s_{t_1} \hat\varphi(\xi)\, d\xi\, ds, \qquad 
	b(t) = -\int\limits^t_{t_1} \hat\varphi(s)\,ds, \qquad c(t) \equiv 1,
\end{equation}
herewith
$t_1$, $t_2 \in [0,T]$ 
are such that 
\begin{equation} \label{14}
	\int\limits^T_0 a(\tau)\, d\tau = 0, 
	\qquad 
	\int\limits^T_0 b(\tau)\, d\tau = 0.
\end{equation}
Due to the conditions on 
$\varphi(t)$, 
the functions 
$a(t)$, $b(t)$ 
are $T$-periodic; i.e.
$$
a(t + T) \equiv a(t), \qquad b(t + T) \equiv b(t).
$$  
	
Using~(\ref{12}) and (\ref{13}),
(\ref{11}) can be rewritten in the following form 
\begin{equation} \label{15}
	\frac{du}{dt} = \mu U(t,\mu)u, 
\end{equation}
$$
U(t,\mu) = U_1 + U_2(t,\mu) + \mu^2U_3(t,\mu),
$$
where
$$
U_1=
\left(
\begin{array}{cc}
	0 & 1
	\\
	- \hat\beta - \frac{1}{T} \int\limits^T_0 \hat\varphi(\tau) a(\tau) d \tau & -\alpha
\end{array}\right),
$$
$$
U_2(t,\mu)= 
\left(
\begin{array}{cc}
	0 & -\mu a(t)
	\\
	-\alpha b(t) - \mu \hat\beta a(t) - \hat\varphi(t)a(t)
	+ \frac{1}{T}\int\limits^T_0 \hat\varphi(\tau) a(\tau) d \tau
	& (\mu a(t) - 1) b(t)
\end{array}\right), 
$$
$$
U_3(t,\mu)=
\left(
\begin{array}{cc}
	0 & \displaystyle{\frac{a^2(t)}{1+\mu a(t)}} 
	\\[7pt]
	0 & -\displaystyle{\frac{a^2(t) b(t)}{1+\mu a(t)}} 
\end{array}\right).
$$
	
As shown in \cite{DeDuMa}, if (\ref{10}) is satisfied, then the spectrum of the matrix 
$U_1$ 
lies in the left half-plane for 
$\alpha > 0$,
from (\ref{14}) it follows that the matrix
$U_2(t, \mu)$
is $T$-periodic and 
\begin{equation} \label{16}
	\int\limits^T_0 U_2(s, \mu)\, ds = 0.
\end{equation}
In this case, as observed in \cite{DeMa01}, \cite{DaKr},  (\ref{15})
belongs to the class of systems whose zero solution is asymptotically stable for sufficiently small 
$\mu$. 

We now indicate a change range of the parameter
$\mu$ 
under which the asymptotic stability holds.
For this purpose we consider the matrix 
$$
{\cal H}(t,\mu) = \frac{1}{\mu} H_1 - H_2(t,\mu),  
$$
where 
$H_1 = H^*_1 > 0$  
is a solution to the Lyapunov matrix equation 
$$
H_1 U_1 + U_1^* H_1 = - I, 
$$
$$
H_2(t,\mu) = H_1 \int\limits^t_0 U_2(s,\mu) ds + \int\limits^t_0 U^*_2(s,\mu) ds\,H_1.
$$ 
It follows from the definition that the matrix 
${\cal H}(t,\mu)$ 
is a Hermitian solution to the boundary value problem of the form~(\ref{6}) 
\begin{equation} \label{17}
	\begin{array}{l}
		\displaystyle
		\frac{d}{dt}{\cal H} + \mu{\cal H}(U_1 + U_2(t,\mu)) + \mu(U_1 + U_2(t,\mu))^*{\cal H} 
		= -C(t,\mu), \quad 0 \le t \le T,\\
		{\cal H}(0,\mu) = {\cal H}(T,\mu), 
	\end{array}
\end{equation}
where 
\begin{equation} \label{18}
	C(t,\mu) = I + \mu H_2(t,\mu) (U_1 + U_2(t,\mu)) + \mu (U_1 + U_2(t,\mu))^* H_2(t,\mu). 
\end{equation}
Taking into account~(\ref{16}), we have 
$$
{\cal H}(0,\mu) = {\cal H}(T,\mu) = \frac{1}{\mu} H_1 > 0.
$$ 
Consequently, by Theorem 2, the zero solution to (\ref{15}) 
with 
$U_3(t,\mu) \equiv 0$ 
is asymptotically stable if    
$C(t,\mu) > 0$, $t \in [0, T]$.
Moreover, it should be noted that the solution to the boundary value problem~(\ref{17}) 
is positive definite for 
$t \in [0, T]$ 
(see \cite{DeMa01, De09}). 
Hence, if we indicate a number 
$\mu_1 > 0$ 
such that 
\begin{equation} \label{19}
	\mu \|H_2(t,\mu) (U_1 + U_2(t,\mu)) + (U_1 + U_2(t,\mu))^* H_2(t,\mu)\| 
	\le \frac{1}{4}, \quad t \in [0, T],
\end{equation}
for 
$\mu \in (0, \mu_1]$,
then 
$$
C(t,\mu) > \frac{3}{4} I, \quad {\cal H}(t,\mu) > 0, \quad t \in [0, T], \quad \mu \in (0, \mu_1]. 
$$
To find 
$\mu_1$,
we estimates firstly the norms of the matrices 
$U_2(t,\mu)$, $U_3(t,\mu)$, $H_2(t,\mu)$. 
	
Introduce the notations 
$$
a = \max\{\alpha, \hat\beta\}, \qquad \varphi_{\max} = \max\limits_{t \in [0,T]} |\hat\varphi(t)|. 
$$
By~(\ref{13}) we have
$$
|a(t)| \le \varphi_{\max} \frac{T^2}{2}, \qquad |b(t)| \le \varphi_{\max} T. 
$$  
Taking into account the form of the matrix 
$U_3(t,\mu)$, 
henceforth we suppose that 
$$
0 < \mu \le \bar\mu = (\varphi_{\max})^{-1} T^{-2}. 
$$ 
Consequently, for 
$\mu \in (0,\, \bar\mu]$,
we have 
\begin{equation} \label{20}
	\|U_3(t,\mu)\| \le \varphi^2_{\max} \frac{T^4}{2} (1 + \varphi_{\max} T),  
\end{equation}
$$
\|U_2(t,\mu)\| \le (1 + a + T)\left(\frac{1}{2} + \varphi_{\max} T\right), 
$$
\begin{equation} \label{21}
	\|H_2(t,\mu)\| \le 2 \|H_1\| T (1 + a + T)\left(\frac{1}{2} + \varphi_{\max} T\right).  
\end{equation}
Therefore, 
\begin{equation} \label{22}
	\|H_2(t,\mu)\|(\|U_1\| + \|U_2(t,\mu)\|) \le L_1,
\end{equation}
where  
$$
L_1 = 2 \|H_1\| T (1 + a + T)\left(\frac{1}{2} + \varphi_{\max} T\right) 
\left(\|U_1\| + (1 + a + T)\left(\frac{1}{2} + \varphi_{\max} T\right)\right).
$$
Define the number 
$$
\mu_1 = \min\left\{\bar\mu, \frac{1}{8L_1}\right\}.  
$$ 
Using the definition~(\ref{18}) of the matrix 
$C(t,\mu)$ 
and~(\ref{22}), we obtain~(\ref{19}) for 
$\mu \in (0, \mu_1]$.
Consequently, 
${\cal H}(t,\mu) > 0$ 
for 
$t \in [0, T]$.
	
Extend the matrix functions  
${\cal H}(t,\mu)$ 
and 
$C(t,\mu)$ 
$T$-periodically on the whole half-axis
$\{t \ge 0\}$, 
keeping the notations.
Let 
$u(t,\mu)$ 
be a solution to~(\ref{15}). Consider the function 
$$
\langle {\cal H}(t,\mu) u(t,\mu), u(t,\mu)\rangle, 
\qquad 
t \ge 0.
$$  
Using the definition of 
${\cal H}(t,\mu)$, 
we obviously have 
$$
\frac{d}{dt} \langle {\cal H}(t,\mu) u(t,\mu), u(t,\mu)\rangle \equiv  
\langle \frac{d}{dt} ({\cal H}(t,\mu)) u(t,\mu), u(t,\mu)\rangle  
$$  
$$
+ \langle {\cal H}(t,\mu) \frac{d}{dt}u(t,\mu), u(t,\mu)\rangle 
+ \langle {\cal H}(t,\mu) u(t,\mu), \frac{d}{dt}u(t,\mu)\rangle \equiv
$$  
$$
- \langle C(t,\mu) u(t,\mu), u(t,\mu)\rangle 
+ \mu^3 \langle ({\cal H}(t,\mu) U_3(t,\mu) + U^*_3(t,\mu) {\cal H}(t,\mu)) u(t,\mu), u(t,\mu)\rangle.
$$  
Introduce the matrix 
$$
{\cal C}(t,\mu) = C(t,\mu) - \mu^3 ({\cal H}(t,\mu) U_3(t,\mu) + U^*_3(t,\mu) {\cal H}(t,\mu)). 
$$ 
Then we can rewrite the previous identity as 
$$
\frac{d}{dt} \langle {\cal H}(t,\mu) u(t,\mu), u(t,\mu)\rangle 
+ \langle {\cal C}(t,\mu) u(t,\mu), u(t,\mu)\rangle \equiv 0. 
$$
	
We now show that there exists 
$\mu_0 \in (0, \mu_1]$ 
such that 
\begin{equation} \label{23}
	{\cal C}(t,\mu) \ge \frac{1}{2} I, \quad t \in [0, T], \quad \mu \in (0, \mu_0].
\end{equation}
Using (\ref{20}) and (\ref{21}), we have 
$$
\mu^3\|{\cal H}(t,\mu) U_3(t,\mu) + U^*_3(t,\mu) {\cal H}(t,\mu)\| 
\le 2\mu^2(\|H_1\| + \mu\|H_2\|)\|U_3\| \le \mu^2 L_2,
$$
for  
$\mu \in (0, \mu_1]$, 
where 
$$
L_2 = \|H_1\| T^4 (\varphi_{\max})^2 (1 + \varphi_{\max} T) 
\left(1 + 2\mu_1 T (1 + a + T)(\frac{1}{2} + \varphi_{\max} T)\right).  
$$
We put 
$$
\mu_0 =\min\{\mu_1, \frac{1}{2\sqrt{L_2}}\}.
$$ 
Then we obtain the inequality 
$$
\mu^3\|{\cal H}(t,\mu) U_3(t,\mu) + U^*_3(t,\mu) {\cal H}(t,\mu)\| \le \frac{1}{4}
$$
for 
$\mu \in (0, \mu_0]$.
Hence, 
$$
{\cal C}(t,\mu) \ge C(t,\mu) - \frac{1}{4}I. 
$$ 
Since 
$C(t,\mu) \ge \frac{3}{4}I$,  
we derive~(\ref{23}). 
	
We now show that the zero solution to (\ref{15}) is asymptotically stable for  
$\mu \in (0, \mu_0]$. 
By (\ref{17}) and (\ref{23}) we have
$$
\frac{d}{dt} \langle {\cal H}(t,\mu) u(t,\mu), u(t,\mu)\rangle 
+ \frac{1}{2} \|u(t,\mu)\|^2 \le 0
$$
for $t \ge 0$.
Taking into account the inequality 
$$
\langle {\cal H}(t,\mu) u(t,\mu), u(t,\mu)\rangle \le \|{\cal H}(t,\mu)\| \|u(t,\mu)\|^2, 
$$ 
we obtain 
$$
\frac{d}{dt} \langle {\cal H}(t,\mu) u(t,\mu), u(t,\mu)\rangle 
+ \frac{1}{2 \|{\cal H}(t,\mu)\|} \langle {\cal H}(t,\mu) u(t,\mu), u(t,\mu)\rangle \le 0. 
$$
Consequently, the following estimate holds   
$$
\langle {\cal H}(t,\mu) u(t,\mu), u(t,\mu)\rangle   
\le \exp\left(-\int\limits^t_0 \frac{1}{2} \|{\cal H}(s,\mu)\|^{-1} ds \right) 
\langle {\cal H}(0,\mu)u(0,\mu), u(0,\mu)\rangle  
$$
for 
$t \ge 0$.
Since 
${\cal H}(t,\mu) > 0$, $t \ge 0$, 
then the zero solution to~(\ref{15}) is asymptotically stable.
	
\vskip10pt

{\bf Lemma.} 
{\it Let the conditions of Theorem 4 be satisfied. Then there exists a number
$\mu_0 > 0$ 
such that the zero solution to~(\ref{9}) is asymptotically stable for every 
$\mu \in (0, \mu_0]$.
}

\vskip10pt
	
The proof of this lemma follows immediately from the fact that we can pass from~(\ref{9})
to~(\ref{11}) which can be transformed into~(\ref{15})
by using the non-degenerate change of variables~(\ref{12}). 
	
\newpage

{\bf 4. Linearized equations with perturbations}

\vskip10pt

We study the case when the coefficients 
$\alpha$, $\hat{\beta}$, $\hat{\varphi}(t)$
of (9) get some perturbations
$$
\alpha \mapsto \alpha + \Delta\alpha, 
\quad  
\hat{\beta} \mapsto \hat{\beta} + \Delta\hat{\beta}, 
\quad 
\hat{\varphi}(t) \mapsto \hat{\varphi}(t) + \Delta\hat{\varphi}(t).
$$
Consider the equation
\begin{equation} \label{24}
\bar{y}'' + (\alpha + \Delta\alpha) \mu \bar{y}' 
+ ((\hat{\beta} + \Delta\hat{\beta}) \mu^2 
+ \mu (\hat{\varphi}(t) + \Delta\hat{\varphi}(t))) \bar{y} = 0,
\qquad
t > 0.
\end{equation}
Pass from the equation to an equivalent system.  
Introduce the vector-function 
$\bar{v}(t,\mu)$ 
such that
$$
\bar{v}(t,\mu) = 
\left(
\begin{array}{c}
	\bar{v}_1(t,\mu) \\ \bar{v}_2(t,\mu)
\end{array}\right)
= \left(\begin{array}{c}
	\bar{y}(t,\mu) \\ \bar{y}'(t,\mu)
\end{array}
\right).
$$
Taking into account (\ref{24}), we obtain 
\begin{equation} \label{25}
\frac{d\bar{v}}{dt} = (A(t,\mu)+\Delta A(t,\mu))\bar{v},
\end{equation}
where
$$
A(t,\mu)+\Delta A(t,\mu) = \left(
\begin{array}{cc}
	0 & 1\\
	-(\hat{\beta} + \Delta\hat{\beta}) \mu^2 - \mu (\hat{\varphi}(t) + \Delta\hat{\varphi}(t)) & -(\alpha + \Delta\alpha)\mu
\end{array}\right).
$$
Henceforth we suppose that the conditions of Theorem 4 are satisfied
and 
$\mu_0 > 0$
is defined in the proof of the lemma.

\vskip10pt
	
{\bf Theorem 5.} 
{\it Let $H(t,\mu)$ be a Hermitian solution to the boundary value problem
\begin{equation} \label{26}
	\begin{array}{l}
		\displaystyle{\frac{d}{dt}}H + HA(t) + A^*(t)H = -I, \quad 0 \le t \le T,
		\\[7pt]
		H(0,\mu) = H(T,\mu)>0,
	\end{array}
\end{equation}
where
$\mu \in (0,\mu_0]$.
If
$\Delta \alpha, \Delta\hat{\beta}, \Delta\hat{\varphi}(t)$
satisfies the inequalities
\begin{equation} \label{27}
	\mu\sup\limits_{t \ge 0} |\Delta\hat{\varphi}(t)|< \frac{1}{4h_{\max}(\mu)},
\end{equation}
\begin{equation} \label{28}
	\mu(|\Delta\hat{\beta}|\mu + |\Delta\alpha|) < \frac{1}{4h_{\max}(\mu)},
\end{equation}
where  
$$
h_{\max}(\mu)=\max\limits_{0 \le t \le T} \|H(t,\mu)\|,
$$
then the zero solution to (\ref{24}) is asymptotically stable.
}

\vskip10pt
	
{\bf Proof.} 
We extend the matrix
$H(t,\mu)$
$T$-periodically on the whole half-axis
$\{t \ge 0\}$, 
keeping the same notation.
Then, for every solution to (\ref{25}), we have
$$
\frac{d}{dt}\langle H \bar{v}, \bar{v} \rangle
= \left\langle \left(\frac{d}{dt} H\right) \bar{v}, \bar{v} \right\rangle
+ \left\langle H \frac{d}{dt} \bar{v}, \bar{v} \right\rangle
+ \left\langle H \bar{v}, \frac{d}{dt} \bar{v} \right\rangle
$$
$$
= \left\langle \left(\frac{d}{dt} H + H A + A^* H\right) \bar{v}, \bar{v} \right\rangle
+ \langle [H \Delta A + \Delta A^* H] \bar{v}, \bar{v} \rangle
$$
$$
= - \|\bar{v}(t)\|^2 + \langle [H \Delta A + \Delta A^* H] \bar{v}, \bar{v} \rangle, \quad t \ge 0.
$$
Hence,
$$
\frac{d}{dt}\langle H \bar{v}(t,\mu), \bar{v}(t,\mu) \rangle
\le -(1-2\|H(t,\mu)\|\|\Delta A(t,\mu)\|)\|\bar{v}(t,\mu)\|^2.
$$
By the conditions of the theorem, we obviously have
\begin{equation} \label{29}
\|\Delta A(t,\mu)\| < \frac{1}{2h_{\max}(\mu)} \le \frac{1}{2\|H(t,\mu)\|}, \quad t \in [0,T].
\end{equation}
Indeed, since
$$
\|\Delta A(t,\mu)\| = \mu\sqrt{(\Delta\hat{\beta}\mu+\Delta\hat{\varphi}(t))^2+(\Delta\alpha)^2}
\le \mu(|\Delta\hat{\beta}|\mu + |\Delta\alpha|) + \mu|\Delta\hat{\varphi}(t)|,
$$
then (\ref{27}) and (\ref{28}) imply (\ref{29}).
Using the inequality
$$
-\|\bar{v}(t,\mu)\| \le -\frac{1}{\|H\|}\langle H(t,\mu) \bar{v}(t,\mu), \bar{v}(t,\mu) \rangle,
$$
we have
$$
\frac{d}{dt}\langle H \bar{v}, \bar{v} \rangle \le 
- \frac{1}{\|H\|}(1-2\|H\|\|\Delta A\|)\langle H \bar{v}, \bar{v} \rangle.
$$
Multiplying this inequality by 
$$
\exp{\left(\int\limits_0^t\left(\frac{1}{\|H(s,\mu)\|} - 2\|\Delta A(s,\mu)\|\right)ds\right)},
$$
we obtain
$$
\frac{d}{dt}\left(\langle H \bar{v}, \bar{v} \rangle \exp{\left(\int\limits_0^t\left(\frac{1}{\|H(s,\mu)\|} 
- 2\|\Delta A(s,\mu)\|\right)ds\right)}\right) \le 0.
$$
Hence,
$$
\langle H \bar{v}, \bar{v} \rangle \exp{\left(\int\limits_0^t\left(\frac{1}{\|H(s,\mu)\|} 
- 2\|\Delta A(s,\mu)\|\right)ds\right)} \le \langle H(0,\mu) \bar{v}(0,\mu), \bar{v}(0,\mu) \rangle.
$$
Consequently,
$$
\|\bar{v}(t,\mu)\|^2 \le \frac{1}{h_{\min}(\mu)}\langle H(0,\mu) \bar{v}(0,\mu), \bar{v}(0,\mu) \rangle
$$
$$
\times
\exp{\left(-\int\limits_0^t\left(\frac{1}{\|H(s,\mu)\|} - 2\|\Delta A(s,\mu)\|\right)ds\right)},
$$
where 
$$
h_{\min}(\mu) \le \min\limits_{t \in [0,T]}\langle H(t,\mu) v, v\rangle, \quad \|v\| = 1.
$$
It follows from this estimate that the zero solution to (\ref{25}) is asympto\-tically stable.
	
The theorem is proven.

\vskip10pt
	
{\bf Corollary.} 
{\it Let the conditions of Theorem 5 be satisfied and let
$H(t, \mu)$
be the $T$-periodic extension 
of the solution to the boundary value problem (\ref{26}) on
$\bar{\mathbb{R}}^+$. 
Then a solution to (\ref{24}) satisfies the estimate
$$
{\|\bar{y}(t,\mu)\|}^2 + \|\bar{y}'_t(t,\mu)\|^2 
\le \frac{\|H(0,\mu)\|}{h_{\min}(\mu)}\left({\|\bar{y}(0,\mu)\|}^2 + \|\bar{y}'_t(0,\mu)\|^2\right)
$$	
$$
\times \exp{\left(-\int\limits_0^t\left(\frac{1}{\|H(s,\mu)\|}-2\|\Delta A(s,\mu)\|\right) d s\right)},
$$
where
{$0<h_{\min}(\mu) \le \min\limits_{t \in [0,T]}h_{\min}(t,\mu)$}
and
$h_{\min}(t,\mu)$ 
is the minimal eigenvalue of the  matrix
$H(t,\mu)$. 
}
	
\vskip20pt

{\bf 5. Stability of stationary solutions to a nonlinear equation with perturbations}

\vskip10pt
	
We now consider a nonlinear differential equation of the form (\ref{1}) 
$$
y'' + \alpha \mu y' + (\beta \mu^2 + \mu \varphi(t)) f(y) = 0, 
$$
where the function 
$\varphi(t)$ 
satisfies (\ref{2}) and the function 
$f(y)$ 
is smooth; moreover,
$$
f(0) = 0, \qquad \frac{d}{dy}f|_{y = 0} < 0.
$$
Hence, 
$y(t) \equiv 0$
is a stationary solution to (\ref{1}). 
We rewrite (\ref{1}) as
$$
y'' + \alpha\mu y' +(\tilde{\beta}\mu^2+\mu\tilde{\varphi}(t))y +(\beta\mu^2+\mu\varphi(t))(f(y) - f'(0)y) = 0,
$$
where
$$
\tilde{\beta} = \beta f'(0), 
\qquad 
\tilde{\varphi}(t) =\varphi(t)f'(0).
$$
We introduce the vector-function 
$x(t)$
as 
$$
x(t) = \left(
\begin{array}{c}
	x_1(t) \\ x_2(t)
\end{array}\right)
= \left(\begin{array}{c}
	y(t) \\ y'(t)
\end{array}
\right).
$$
Then we obtain the following system equivalent to (\ref{1})
$$
\frac{d}{dt}x = \tilde{A}(t,\mu)x + F(t, \mu, x),
$$
where
$$
\tilde{A}(t,\mu) = \left(
\begin{array}{cc}
	0 & 1\\
	-\tilde{\beta} \mu^2 - \mu \tilde{\varphi}(t) & -\alpha\mu
\end{array}\right),
$$
$$
F(t, \mu, x) = \left(
\begin{array}{c}
	0 \\ -(\beta \mu^2 + \mu \varphi(t))(f(x_1)-f'(0)x_1)
\end{array}\right).
$$	
We now consider the nonlinear equation with perturbations of the coefficients
$ \alpha, \beta, \varphi(t)$
\begin{equation} \label{30}
	\bar{y}'' + (\alpha + \Delta\alpha) \mu \bar{y}' + ((\beta + \Delta\beta) \mu^2 + \mu (\varphi(t) + \Delta\varphi(t))) f(\bar{y}) = 0.
\end{equation}
Pass from the equation to an equivalent system in a similar way by introducing the vector-function
$$
\bar{x}(t) = \left(
\begin{array}{c}
	\bar{x}_1(t) \\ \bar{x}_2(t)
\end{array}\right)
= \left(
\begin{array}{c}
	\bar{y}(t) \\ \bar{y}'(t)
\end{array}
\right).
$$
We obtain the following system
\begin{equation} \label{31}
	\frac{d}{dt}\bar{x} = (\tilde{A}+\Delta \tilde{A})(t,\mu)\bar{x} + F(t, \mu, \bar{x}),
\end{equation}
where 
$$
\tilde{A}(t,\mu) + \Delta \tilde{A}(t,\mu) = \left(
\begin{array}{cc}
	0 & 1\\
	-(\tilde{\beta} + \Delta \tilde{\beta}) \mu^2 - \mu (\tilde{\varphi}(t) + \Delta \tilde{\varphi}(t))  & -(\alpha+ \Delta \alpha)\mu
\end{array}\right),
$$
$$
F(t, \mu, \bar{x}) = \left(
\begin{array}{c}
	0 \\ -[(\beta+\Delta\beta) \mu^2 + \mu (\varphi(t)+\Delta\varphi(t))](f(\bar{x}_1)-f'(0)\bar{x}_1)
\end{array}\right).
$$
To study the asymptotic stability of the zero solution to (\ref{30}),
we consider the Cauchy problem
\begin{equation} \label{32}
\begin{array}{l}
	\tilde{y}'' + (\alpha + \Delta\alpha) \mu \tilde{y}' + ((\beta + \Delta\beta) \mu^2 + \mu (\varphi(t) + \Delta\varphi(t))) f(\tilde{y}) = 0,
	\\
	\tilde{y}|_{t=0}=y^0,
	\\
	\tilde{y}'|_{t=0}=y^1.
\end{array}
\end{equation}
Henceforth we suppose that the conditions of Theorem 5 are satisfied,
$$
\mu \in (0,\mu_0], \quad |f(\xi)-f'(0)\xi| \le p{|\xi|}^2, \quad \xi
\in \mathbb{R},
$$
$$
q(\mu) = \sup\limits_{t \ge 0}[(|\beta +\Delta\beta|\mu^2 + |\varphi(t)+\Delta\varphi(t)|\mu)p].
$$
Let
$H(t,\mu)$
be a solution to the boundary value problem
$$
\begin{array}{l}
	\displaystyle{\frac{d}{dt}}H + H\tilde{A}(t,\mu) + \tilde{A}^*(t,\mu)H = -I, \quad 0 \le t \le T,
	\\[7pt]
	H(0, \mu) = H(T, \mu)>0.
\end{array}
$$
We extend the matrix
$H(t,\mu)$
$T$-periodically on the whole half-axis
$\{t \ge 0\}$, 
keeping the same notation.
For every solution to (\ref{31}), we have
$$
\frac{d}{dt}\langle H \bar{x},\bar{x} \rangle
= \left\langle \left(\frac{d}{dt} H\right) \bar{x}, \bar{x} \right\rangle
+ \left\langle H \frac{d}{dt}\bar{x},\bar{x} \right\rangle
+ \left\langle H \bar{x}, \frac{d}{dt} \bar{x} \right\rangle
$$
$$
= \left\langle \left(\frac{d}{dt} H + H \tilde{A} + \tilde{A}^* H\right) \bar{x}, \bar{x} \right\rangle
+ 2\langle H \Delta\tilde{A}\bar{x}, \bar{x} \rangle + 2\langle H F(t,\mu, \bar{x}), \bar{x}\rangle, \quad t \ge 0.
$$
Hence,
$$
\frac{d}{dt}\langle H \bar{x},\bar{x} \rangle \le - \|\bar{x}\|^2 + 2\|H\|\| \Delta\tilde{A}\|\|\bar{x}\|^2 + 2\|H \|\|F(t,\mu,\bar{x})\|\|\bar{x}\|
$$
$$
\le - (1-2\|H\|\| \Delta\tilde{A}\|)\|\bar{x}\|^2 + 2q(\mu)h_{\max}(\mu)\|\bar{x}\|^3.
$$
Untroducing the notations
\begin{equation} \label{33}
	\tilde{q}(\mu)= 2q(\mu)h_{\max}(\mu),  
\quad 
	\varepsilon(t,\mu)=1-2\|H(t,\mu)\|\| \Delta\tilde{A}(t,\mu)\|, 
\end{equation}
we obtain
$$
\frac{d}{dt}\langle H \bar{x},\bar{x} \rangle \le - \varepsilon(t,\mu)\|\bar{x}\|^2 + \tilde{q}(\mu)\|\bar{x}\|^3.
$$
By the inequalities
$$
h_{\min}(\mu)\|\bar{x}\|^2 \le \langle H \bar{x},\bar{x} \rangle \le \|H(t,\mu)\|\|\bar{x}\|^2, 
\quad 
\varepsilon(t,\mu) > 0,
$$
we have
$$
\frac{d}{dt}\langle H \bar{x},\bar{x} \rangle 
\le - \frac{\varepsilon(t,\mu)}{\|H(t,\mu)\|} \langle H \bar{x},\bar{x} \rangle 
+ \frac{\tilde{q}(\mu)}{h_{\min}^{3/2}(\mu)}\langle H \bar{x},\bar{x} \rangle^{3/2}.
$$
Divide both parts by 
$\langle H \bar{x},\bar{x} \rangle^{3/2}$
for $\bar{x} \neq 0$ 
Then,
$$
\frac{1}{\langle H \bar{x},\bar{x} \rangle^{3/2}}\frac{d}{dt}\langle H \bar{x},\bar{x} \rangle 
+ \frac{\varepsilon(t,\mu)}{\|H(t,\mu)\|}\langle H \bar{x},\bar{x} \rangle^{-1/2} 
\le \frac{\tilde{q}(\mu)}{h_{\min}^{3/2}(\mu)}.
$$
Multiplying this inequality by
$$ 
\exp{\left(-\int\limits_0^t\frac{\varepsilon(\tau,\mu)}{2\|H(\tau,\mu)\|}d\tau\right)}, 
$$ 
we obtain
$$
\frac{d}{dt}
\left[\exp{\left(-\int\limits_0^t
\frac{\varepsilon(\tau,\mu)}{2\|H(\tau,\mu)\|}d\tau\right)}
\langle H \bar{x},\bar{x} \rangle^{-1/2}\right]
\ge - \frac{\tilde{q}(\mu)}{2 h_{\min}^{3/2}(\mu)}
\exp{\left(-\int\limits_0^t\frac{\varepsilon(\tau,\mu)}{2\|H(\tau,\mu)\|}d\tau\right)}.
$$
Consequently, 
$$
\exp{\left(-\int\limits_0^t
\frac{\varepsilon(\tau,\mu)}{2\|H(\tau,\mu)\|}d\tau\right)}
\langle H \bar{x},\bar{x} \rangle^{-1/2} - \langle H \bar{x}(0),\bar{x}(0) \rangle^{-1/2}
$$
$$
\ge - \frac{\tilde{q}(\mu)}{2 h_{\min}^{3/2}(\mu)}
\int\limits_0^t\exp\left(-\int\limits_0^\tau
\frac{\varepsilon(\xi,\mu)}{2\|H(\xi,\mu)\|}d\xi\right)d\tau.
$$
Let 
$\psi(t,\mu) = \langle H(t,\mu) \bar{x}(t,\mu),\bar{x}(t,\mu) \rangle$. 
Then,
$$
\exp{\left(-\int\limits_0^t\frac{\varepsilon(\tau,\mu)}{2\|H(\tau,\mu)\|}d\tau\right)}
\psi^{-1/2}(t,\mu)
$$
$$
\ge \psi^{-1/2}(0,\mu) - \frac{\tilde{q}(\mu)}{2 h_{\min}^{3/2}(\mu)}
\int\limits_0^t\exp{\left(-\int\limits_0^\tau\frac{\varepsilon(\xi,\mu)}{2\|H(\xi,\mu)\|}d\xi\right)d\tau}.
$$
Multiplying this inequality by 
$\psi^{1/2}(0,\mu)\psi^{1/2}(t,\mu)$,
we derive
$$
\exp{\left(-\int\limits_0^t\frac{\varepsilon(\tau,\mu)}{2\|H(\tau,\mu)\|}d\tau\right)}
\psi^{1/2}(0,\mu)
$$
$$
\ge \left[1- \frac{\tilde{q}(\mu)}{2 h_{\min}^{3/2}(\mu)}\psi^{1/2}(0,\mu)\int\limits_0^t\exp{\left(-\int\limits_0^\tau\frac{\varepsilon(\xi,\mu)}{2\|H(\xi,\mu)\|}d\xi\right)d\tau}\right]\psi^{1/2}(t,\mu).
$$
	
Consider the following integral 
$$ 
I(t,\mu) =  \int\limits_0^t\exp{\left(-\int\limits_0^\tau
\frac{\varepsilon(\xi,\mu)}{2\|H(\xi,\mu)\|}d\xi\right)d\tau}. 
$$
Proving Theorem 5, by (\ref{27}) and (\ref{28}) we established (\ref{29}) 
for the norm of the matrix
$\Delta\tilde{A}(\xi,\mu)$:
$$
\|\Delta\tilde{A}(\xi,\mu)\| < \frac{1}{2h_{\max}(\mu)}.
$$
Further we need a stronger estimate.
Therefore we suppose that, instead of (\ref{27}) and (\ref{28}),
the following conditions are fulfilled
\begin{equation} \label{34}
	\mu\sup\limits_{t \ge 0} |\Delta\tilde{\varphi}(t)|< \frac{1}{8h_{\max}(\mu)},
\end{equation}
\begin{equation} \label{35}
	\mu(|\Delta\tilde{\beta}|\mu + |\Delta\alpha|) < \frac{1}{8h_{\max}(\mu)}.
\end{equation}
Then,
$$
\|\Delta\tilde{A}(t,\mu)\| = \mu\sqrt{(\Delta\tilde{\beta}\mu + \Delta\tilde{\varphi}(t))^2 + (\Delta\alpha)^2}
\le \mu(\Delta\tilde{\beta}\mu+\Delta\alpha) + \mu|\Delta\tilde{\varphi}(t)| \le \frac{1}{4h_{\max}(\mu)}.
$$
Using this inequality, we have
$$
\varepsilon(\xi,\mu) = 1-2\|H(\xi,\mu)\|\| \Delta\tilde{A}(\xi,\mu)\| \ge 1 - 2h_{\max}(\mu)\|\Delta\tilde{A}(\xi,\mu)\| \ge \frac{1}{2}.
$$
Consequently, we obtain 
$$
|I(t,\mu)| \le \int\limits_0^t\exp\left(-\frac{\tau}{4h_{\max}(\mu)}\right)d\tau \le 4h_{\max}(\mu).
$$

Hence,
$$
\exp{\left(-\int\limits_0^t\frac{\varepsilon(\tau,\mu)}{2\|H(\tau,\mu)\|}d\tau\right)}
\psi^{1/2}(0,\mu)
\ge \left[1- \frac{\tilde{q}(\mu)}{ h_{\min}^{3/2}(\mu)}\psi^{1/2}(0,\mu)2h_{\max}(\mu)\right]
\psi^{1/2}(t).
$$
Suppose that 
$$
\langle H(0, \mu) \bar{x}(0,\mu),\bar{x}(0,\mu) \rangle \le \frac{h_{\min}^3(\mu)}{64q^2h_{\max}^4(\mu)}.
$$
Then,
$$
\|\bar{x}(t,\mu)\|^2 \le \frac{1}{h_{\min}(\mu)}\langle H(t,\mu) \bar{x}(t,\mu),\bar{x}(t,\mu) \rangle
$$
$$
\le \frac{4}{h_{\min}(\mu)}\exp{\left(-\int\limits_0^t \frac{1}{\|2H(\tau,\mu)\|}d\tau\right)}\langle H(0,\mu) \bar{x}(0),\bar{x}(0) \rangle.
$$
	
Taking into account the above reasoning, we conclude the following results.

\vskip10pt
	
{\bf Theorem 6.} 	
{\it Let the conditions of Theorem 5 be satisfied and let the pertur\-bations 
{\large$\Delta\alpha, \Delta\beta, \Delta\varphi(t)$}
satisfy (\ref{34}) and (\ref{35}).
If the components of the initial vector
$y^0, y^1$
are such that the inequality holds
\begin{equation} \label{36}
	\left\langle H(0,\mu)\left(
	\begin{array}{c}
		y^0 \\ y^1
	\end{array}\right), \left(
	\begin{array}{c}
		y^0 \\ y^1
	\end{array}\right) \right\rangle \le \frac{h^3_{\min}(\mu)}{64h_{\max}^4(\mu)q^2(\mu)},
\end{equation}
$$
q(\mu) = \sup\limits_{t \ge 0}[(|\beta +\Delta\beta|\mu^2 + |\varphi(t)+\Delta\varphi(t)|\mu)p],
$$
then the solution to the Cauchy problem (\ref{32}) is uniquely defined for 
$t \ge 0$.
}

\vskip10pt
	
{\bf Theorem 7.} 	
{\it Let the conditions of Theorem 6 be satisfied.
Then the solution to the Cauchy problem (\ref{32}) satisfies the estimate
$$
{\|\tilde{y}(t,\mu)\|}^2 + \|\tilde{y}'_t(t,\mu)\|^2
$$
\begin{equation} \label{37}
	\le\frac{4}{h_{\min}(\mu)} \exp\left({-\int\limits^t_0\frac{1}{2\|H(\tau,\mu)\|} d\tau}\right)  \left\langle H(0,\mu)\left(
	\begin{array}{c}
		y^0 \\ y^1
	\end{array}\right), \left(
	\begin{array}{c}
		y^0 \\ y^1
	\end{array}\right) \right\rangle.
\end{equation}
}	

\vskip10pt

{\bf Corollary.}
{\it If the conditions of Theorem 6 are satisfied,
then the zero solution to (\ref{30}) is asymptotically stable.
}

\vskip10pt
	
{\bf Remark.}
(\ref{36}) gives an estimate for the attraction set of the zero solution to (\ref{30})
and (\ref{37}) gives an estimate for the stabilization rate of the solution as
$t \to \infty$.

\vskip10pt

{\bf Theorem 8.} 
{\it Let the conditions of Theorem 6 be satisfied,
$$
|f(\xi)-f'(0)\xi| \le p{|\xi|}^2, 
\qquad |\xi| \le \rho.
$$
If the components of the initial vector
$y^0, y^1$
are such that the inequalities hold
\begin{equation} \label{38}
	\left\langle H(0,\mu)\left(
	\begin{array}{c}
		y^0 \\ y^1
	\end{array}\right), \left(
	\begin{array}{c}
		y^0 \\ y^1
	\end{array}\right) \right\rangle \le \frac{h^3_{\min}(\mu)}{64h_{\max}^4(\mu)q^2(\mu)},
\end{equation}
\begin{equation} \label{39}
\sqrt{(y^0)^2 + (y^1)^2} \le \frac{\rho h_{\min}(\mu)}{4h_{\max}(\mu)},
\end{equation}
then the solution to the Cauchy problem (\ref{32}) is uniquely defined for 
$t \ge 0$
and it satisfies (\ref{37}).
}

\vskip10pt
	
{\bf Corollary.}
{\it If the conditions of Theorem 8 are satisfied,
then the zero solution to (\ref{30}) is asymptotically stable.
}

\vskip10pt
	
{\bf Remark.} 
(\ref{38}) and (\ref{39}) give estimates for the attraction set of the zero solution 
to (\ref{30}).

\vskip10pt

{\bf Acknowledgments.}
The study was carried out within the framework of the state contract
of the Sobolev Institute of Mathematics (project no.~FWNF-2022-0008).

\end{document}